\input amssym.def
\input amssym
\magnification=1200
\parindent0pt
\hsize=16 true cm
\baselineskip=13  pt plus .2pt

$ $

\def\F{{\cal F}}
\def\R{\Bbb R}
\def\Z{\Bbb Z}

\def\S{\Bbb S}

\centerline {\bf  A note on actions of the  symplectic group ${\rm Sp}(2g,\Bbb
Z)$ on homology spheres}

\bigskip \bigskip

\centerline {Bruno P. Zimmermann}

\bigskip

\centerline {Universit\`a degli Studi di Trieste}

\centerline {Dipartimento di Matematica e Informatica}

\centerline {34127 Trieste, Italy}

\centerline {zimmer@units.it}

\vskip 1cm

Abstract. {\sl  The symplectic group ${\rm Sp}(2g,\Bbb Z)$ is a subgroup of the
linear  group  ${\rm SL}(2g,\Bbb Z)$ and  admits a faithful action on the sphere
$S^{2g-1}$, induced from its linear action on Euclidean space $\Bbb R^{2g}$.
Generalizing corresponding results for linear groups, we show that, if
$m < 2g-1$ and $g> 2$, any continuous action of ${\rm Sp}(2g,\Bbb Z)$ on a
homology $m$-sphere, and in particular on  $S^m$, is trivial. }

\bigskip \bigskip

{\bf 1. Introduction}

\medskip

The linear group ${\rm SL}(n,\Bbb Z)$ admits a faithful, linear action on
Euclidean space $\R^n$, on the sphere $S^{n-1}$ and on the torus $T^n = \Bbb
R^n/\Bbb Z^n$. On the other hand, by [We] any smooth action of ${\rm SL}(n,\Bbb
Z)$ on the torus $T^m$ is trivial if $m < n$ and $n > 2$, and by results of
Parwani [Pa] and the author [Z1], any smooth action of the linear group ${\rm
SL}(n,\Bbb Z)$ on a mod 2 homology $m$-sphere is trivial if $m < n-1$ and $n >
2$. We note that the results in [Pa] are stated for arbitrary continuous
actions but it is noticed in [BV, Remarks 4.16 and 4.17] that some of the
arguments involving Smith fixed point theory are not correct in this more
general setting; in [BV], the result is obtained for continuous actions as a
consequence of a corresponding result for actions of the automorphism group of
a free group of rank $n$ on homology spheres.

\medskip

Our main  result is the following analogue for actions of  the symplectic group
${\rm Sp}(2g,\Bbb Z)$ (which admits a faithful, linear action on $S^{2g-1}$).

\bigskip

{\bf Theorem}. {\sl  For $m < 2g-1$ and $g \ge 3$, any continuous action of the
symplectic group ${\rm Sp}(2g,\Bbb Z)$ on a homology $m$-sphere, and in
particular on $S^m$, is trivial (for $g = 2$, it factors through the action of a
finite group). }

\bigskip

In fact, we will prove the Theorem for mod 3 homology $m$-spheres, i.e. for
manifolds with the mod 3 homology of the $m$-sphere (homology with coefficients
in the integers mod 3).

\medskip

It should be noted that the proofs of most of these results depend strongly on
the existence of torsion, that is of certain types of finite subgroups in the
groups considered (3-torsion in the case of the Theorem), and that the real
challenge is to prove such results for arbitrary, possibly torsion-free
subgroups of finite index. For the case of $S^1$, Witte [Wi] has shown that
every continuous action of a subgroup of finite index in ${\rm SL}(n,\Bbb Z)$
on the circle $S^1$ factors through a finite group action, for $n > 2$. In the
context of the Zimmer program for actions of irreducible lattices in semisimple
Lie groups  of $\Bbb R$-rank at least two, it is conjectured by Farb and Shalen
([FS]) that any smooth action of a finite-index subgroup of ${\rm SL}(n,\Bbb
Z)$ on a compact $m$-manifold factors through the action of a finite group if
$m < n-1$ and $n > 2$.

\medskip

Abelianization of the fundamental group $\pi_1(\F_g)$ of a closed, orientable
surface $\F_g$ of genus $g$ induces a surjection
$${\rm MC}(\F_g) \to {\rm Sp}(2g,\Bbb Z)$$ of the mapping class group ${\rm
MC}(\F_g) \cong {\rm Out}_+(\pi_1(\F_g))$ of isotopy classes of
orientation-preserving homeomorphisms of $\F_g$ onto the symplectic group ${\rm
Sp}(2g,\Bbb Z)$, and the following question naturally arises: What is the
minimal dimension of a non-trivial (or infinite image, or faithful) action of
the mapping class group ${\rm MC}(\F_g)$ on a sphere or a homology sphere?  See
[BV] for the case of the outer automorphism group  ${\rm Out}_+(F_n)$ of a free
group $F_n$ of rank $n$. The case of mapping class groups appears more
difficult, due to the irregular distribution of torsion depending on the genus;
we present the following partial result for genus three. Note that, via the
surjection ${\rm MC}(\F_g) \to {\rm Sp}(2g,\Bbb Z)$, the mapping class group
${\rm MC}(\F_g)$ admits a non-trivial, linear action on $S^{2g-1} \subset
\R^{2g}$.

\bigskip

{\bf Proposition}. {\sl  Any smooth action of the mapping class group ${\rm
MC}(\F_3)$ on a homology sphere of dimension less than five is trivial.}

\bigskip

But even in this case the minimal dimension of a faithful action, or of an
action with infinite image, remains open; note that ${\rm MC}(\F_3)$ acts
faithfully on the Teichm\"uller space in genus three, homeomorphic to
$\R^{12}$, and on its boundary homeomorphic to $S^{11}$.

\bigskip

{\bf 2. Proof of the Theorem}

\medskip

Our main reference for the symplectic group ${\rm Sp}(2g,\Bbb Z)$ is the book of
Newman [N].

\medskip

Fixing a basis $a_1, b_1, \ldots, a_g, b_g$ of the free abelian group $\Z^{2g}$,
one considers the anti-symmetric or symplectic bilinear form on $\Z^{2g}$
defined by
$$a_i\times b_i = 1, \;\;  a_i \times b_j = 0 \;\; {\rm for} \;\; i \ne j,
\;\;  a_i \times a_j = 0,$$ that is with matrix
$$J =  \pmatrix {
                 0  &  1   \cr
                 -1  &  0    \cr } \oplus \ldots  \oplus
                \pmatrix {
                 0  &   1    \cr
                 -1  &  0    \cr }$$ with respect to the chosen basis (so $J$ is
a $(2g \times 2g)$ matrix consisting of $g$ blocks).

\medskip

The symplectic group ${\rm Sp}(2g,\Bbb Z)$ is the group of automorphisms of
$\Z^{2g}$ which preserve this bilinear form, or equivalently the subgroup of
integer
$2g \times 2g$ matrices $A$ of the linear group ${\rm SL}(2g,\Bbb Z)$ such that
$AJA^T=J$ (where $A^T$ denotes the transposed matrix); note that
${\rm Sp}(2,\Bbb Z) = {\rm SL}(2,\Bbb Z)$.

\medskip

For a fixed $i$, we consider the symplectic automorphism $\phi_i$ of order three
of
$\Z^{2g}$  which fixes all $a_j, b_j$ with $j \ne i$, and with
$\phi_i(a_i) = -a_i - b_i$,   $\phi_i(b_i) = a_i$, that is with matrix
$$\pmatrix {
                 -1  &  1   \cr
                 -1  &  0    \cr }$$ of order three on the subgroup $\Z^2 =
\langle a_i, b_i \rangle$ of $\Z^{2g}$. We denote by  $\Z_3^g$ the elementary
abelian subgroup of ${\rm Sp}(2g,\Bbb Z)$  generated by the symplectic
automorphisms $\phi_1, \ldots, \phi_g$.

\medskip

Let $g \ge 3$. We consider an action of ${\rm Sp}(2g,\Bbb Z)$ on a mod 3
homology sphere $M$ of dimension less than $2g - 1$, and have to show that the
action is trivial.  We denote by $U$ the normal subgroup of ${\rm Sp}(2g,\Bbb
Z)$ consisting of all elements which act trivially on $M$. Crucial for the
proof of the Theorem is the following well-known result from Smith fixed point
theory ([S]; see also [BV] for a useful discussion of the concept of a
generalized manifold suited for the proofs).

\bigskip

{\bf Proposition 1.}  {\sl  For a prime $p$, the minimal dimension of a
faithful, continuous action of the elementary abelian group $\Bbb Z_p^k$  on a
mod $p$ homology sphere is $2k-1$ if $p$ is odd, and $k-1$ if $p = 2$.}

\bigskip

By  Proposition 1, some element $u$ of order three in $\Z_3^g$ acts trivially
on $M$, that is lies in the kernel $U$ of the action of ${\rm Sp}(2g,\Bbb Z)$
on $M$; note that $u$ is non-central in ${\rm Sp}(2g,\Bbb Z)$. By the
congruence subgroup property for the symplectic groups ([Me], and [BMS] for a
more general version), for $g > 1$ any non-central, normal subgroup $U$ of the
symplectic group ${\rm Sp}(2g,\Bbb Z)$ has finite index and contains a
congruence subgroup $C(k)$, for some positive integer $k$; here $C(k)$ denotes
the kernel of the canonical map ${\rm Sp}(2g,\Bbb Z) \to {\rm Sp}(2g,\Z/k\Z)$
which is surjective by [N, Theorem VII.21], so we have an exact sequence
$$1 \to  C(k)  \to  {\rm Sp}(2g,\Bbb Z) \to {\rm Sp}(2g,\Z/k\Z) \to 1.$$
Hence the action of ${\rm Sp}(2g,\Bbb Z)$ on $M$
factors through the action of a finite group $G$, with $G \cong {\rm Sp}(2g,\Bbb
Z)/U$,  and  there is a surjection from the finite group ${\rm Sp}(2g,\Bbb
Z)/C(k) \cong {\rm Sp}(2g,\Bbb Z/k\Bbb Z)$ to ${\rm Sp}(2g,\Bbb Z)/U \cong G$
which we denote by $$\Phi: {\rm Sp}(2g,\Z/k\Z) \to G.$$

\bigskip

Remark. In order to see that the action of ${\rm Sp}(2g,\Bbb Z)$ on $M$ factors
through the action of a finite group $G$, one may apply also the Margulis
finiteness theorem which states that an irreducible lattice in a semisimple Lie
group of real rank at least two is almost simple, that is any normal subgroup of
the lattice is either of finite index, or contained in the center of the
semisimple Lie group (see [Ma] or [Z, Theorem 8.1.2]).  Note that the Margulis
finiteness theorem applies to the irreducible lattice  ${\rm Sp}(2g,\Bbb Z)$ in
the semisimple Lie group ${\rm Sp}(2g,\Bbb R)$, for $g \ge 2$.

\bigskip

We have to show that the finite quotient $G$ of ${\rm Sp}(2g,\Bbb Z)$ is
trivial; assume, by contradiction, that it is not. Then, since ${\rm
Sp}(2g,\Bbb Z)$ is perfect for $g \ge 3$, also $G$ is perfect and hence
non-solvable (more generally, by [Po], ${\rm MC}(\F_g)$ is perfect for $g \ge
3$ whereas the abelianization of ${\rm MC}(\F_2)$ is cyclic of order 10).

\medskip

Considering the prime decomposition  $k = p_1^{r_1} \ldots p_s^{r_s}$ of $k$,
one has
$${\rm Sp}(2g,\Bbb Z/k\Bbb Z) \; \cong \; {\rm Sp}(2g,\Bbb Z/p_1^{r_1}\Z)
\times \ldots \times {\rm Sp}(2g,\Bbb Z/p_s^{r_s}\Z)$$ (see [N, Theorem VII.26]).
The restriction of $\Phi: {\rm Sp}(2g,\Z/k\Z) \to G$ to some factor
${\rm Sp}(2g,\Bbb Z/p_i^{r_i}\Z) = {\rm Sp}(2g,\Bbb Z/p^r\Z)$ has to be
non-trivial and induces a surjection $$\Phi_0: {\rm Sp}(2g,\Z/p^r\Z) \to G_0$$
onto a perfect, non-solvable subgroup $G_0$ of $G$; we denote by $U_0$ the
kernel of $\Phi_0$ (the elements of ${\rm Sp}(2g,\Z/p^r\Z)$ acting trivially on
$M$).

\medskip

Let $K$ denote the kernel of the canonical surjection ${\rm Sp}(2g,\Z/p^r\Z) \to
{\rm Sp}(2g,\Z/p\Z)$, so $K$ consists of all matrices in ${\rm Sp}(2g,\Z/p^r\Z)$
which are congruent to the identity matrix $I = I_{2g}$ when entries are taken
mod
$p$. By performing the binomial expansion of $(I + pA)^{p^{r-1}}$ one checks
that $K$ is a $p$-group, in particular $K$ is solvable.

\medskip

Let $K_0$ denote the kernel of the surjection from ${\rm Sp}(2g,\Z/p^r\Z)$ to
the central quotient ${\rm PSp}(2g,\Z/p\Z)$ of ${\rm Sp}(2g,\Z/p\Z)$; also
$K_0$ is solvable and, since $g \ge 3$, ${\rm PSp}(2g,\Z/p\Z)$ is a non-abelian
simple group (${\rm PSp}(4,\Z/2\Z)$ is isomorphic to the symmetric group
$\S_6$).

\medskip

Note that the element $u$ of order three in ${\rm Sp}(2g,\Z)$ considered above
injects into the successive quotients ${\rm Sp}(2g,\Z/k\Z)$, ${\rm Sp}(2g,\Bbb
Z/p^r\Z)$, ${\rm Sp}(2g,\Z/p\Z)$ and ${\rm PSp}(2g,\Z/p\Z)$, hence the normal
subgroup $U_0$ of ${\rm Sp}(2g,\Z/p^r\Z)$ surjects onto the finite simple group
${\rm PSp}(2g,\Z/p\Z)$.

\medskip

We consider the two exact sequences

$$ 1 \to K_0 \to {\rm Sp}(2g,\Z/p^r\Z) \to {\rm PSp}(2g,\Z/p\Z) \to 1,$$

$$ 1 \to U_0 \cap K_0 \to U_0 \to {\rm PSp}(2g,\Z/p\Z) \to 1.$$

Quotienting the first by the second, we conclude that  ${\rm
Sp}(2g,\Z/p^r\Z)/U_0
\cong G_0$ is isomorphic to the solvable group $K_0/(U_0 \cap K_0)$. This is a
contradiction, and hence  $G$ has to be trivial.

\medskip

This completes the proof of the Theorem.

\bigskip

{\bf 2. Proof of the Proposition}

\medskip

The mapping class group  ${\rm MC}(\F_3)$ has a finite  subgroup isomorphic to
the linear fractional group ${\rm PSL}(2,\Z/7\Z)$ of order 168, induced from
the Hurwitz action of the smallest Hurwitz group ${\rm PSL}(2,\Z/7\Z)$ on
Klein's quartic of genus three (i.e., of maximal possible order $84(g-1)$). By
[Z3, Proposition 1] and [MZ, Proposition 2.1], the simple group ${\rm
PSL}(2,\Z/7\Z)$ does not admit a non-trivial, smooth action on a homology
sphere of dimension less than five (see also [Z4]). On the other hand, by [Z2,
Lemma 1] any non-trivial homomorphism from ${\rm MC}(\F_3)$ to an arbitrary
group has to inject  ${\rm PSL}(2,\Z/7\Z)$, hence any action of  ${\rm
MC}(\F_3)$ on a homology sphere of dimension less than five is trivial.

\bigskip
\bigskip

\centerline {\bf References}

\bigskip

\item {[BMS]} H.Bass, J.Milnor, J.P.Serre, {\it The congruence subgroup property
for
$SL_n$ ($n \ge 3$) and  $SP_{2n}$ ($n \ge 2$).}    Inst. Hautes Etudes Sci. Publ.
Math. 33, 59-137 (1967)

\item {[BV]} M.R.Bridson, K.Vogtmann, {\it  Actions of automorphism groups of
free groups on homology spheres and acyclic manifolds.}  Electronic version in
arXiv:0803.2062

\item {[FS]} B.Farb, P.Shalen, {\it Real-analytic actions of lattices.} Invent.
math. 135, 273-296  (1999)

\item {[Ma]} G.Margulis, {\it Discrete Subgroups of Semisimple Lie Groups.}
Ergeb. Math. Grenzgebiete 17,  Springer-Verlag 1991

\item {[Me]} J.L.Mennicke, {\it Zur Theorie der Siegelschen Modulgruppe.} Math.
Ann.  159, 115-129  (1965)

\item {[MZ]} M.Mecchia, B.Zimmermann, {\it On finite simple and nonsolvable
groups acting on homology 4-spheres.} Top. Appl. 153,  2933-2942  (2006)

\item {[N]} M.Newman, {\it Integral Matrices.} Pure and Applied Mathematics
Vol.45, Academic Press 1972

\item {[Pa]} K.Parwani, {\it Actions of SL$(n,\Bbb Z)$ on homology spheres.}
Geom. Ded. 112, 215-223 (2005)

\item {[Po]} J.Powell, {\it Two theorems on the mapping class group of a
surface.}  Proc. Amer. Math. Soc. 68, 347-350  (1978)

\item {[S]} P.A.Smith, {\it Permutable periodic transformations.} Proc. Nat.
Acad. Sci. U.S.A. 30, 105 - 108 (1944)

\item {[We]} S.Weinberger, {\it SL$(n,\Bbb Z)$ cannot act on small tori.} AMS/IP
Studies in Advanced Mathematics Volume 2 (Part 1), 406-408 (1997)

\item {[Wi]} D.Witte, {\it Arithmetic groups of higher $\Bbb Q$-rank cannot act
on 1-manifolds.}  Proc. Amer. Math. Soc. 122, 333-340 (1994)

\item {[Z]} R.Zimmer, {\it Ergodic Theory and Semisimple Groups.} Monographs in
Mathematics 81,  Birkh\"auser 1984

\item {[Z1]} B.Zimmermann, {\it ${\rm SL}(n,\Bbb Z)$ cannot act on small
spheres.} Top. Appl. 156,  1167-1169  (2009) (arXiv:math/0604239)

\item {[Z2]} B.Zimmermann, {\it  A note on minimal finite quotients of mapping
class groups.}   Electronic version in  arXiv:0803.3144

\item {[Z3]} B.Zimmermann, {\it On finite simple groups acting on homology
3-spheres.}  Top. Appl. 125, 199-202 (2002)

\item {[Z4]} B.Zimmermann, {\it On the minimal dimension of a homology sphere on
which a finite group acts.}  Math. Proc. Camb. Phil. Soc. 144, 397-401 (2008)

\bye